\documentclass[a4paper]{article}
\usepackage{amsthm,amsmath,amssymb}
\usepackage{color}

\newtheorem{teor}{Theorem}[section]
\newtheorem{defin}[teor]{Definition}
\newtheorem{lemm}[teor]{Lemma}
\newtheorem{osse}[teor]{Remark}
\newtheorem{assu}[teor]{Assumption}
\newtheorem{prop}[teor]{Proposition}
\newtheorem{defi}[teor]{Definition}
\newtheorem{coro}[teor]{Corollary}
\newtheorem{prob}[teor]{Problem}

\newcommand{\bele}{\begin{lemm}\begin{sl}}
\newcommand{\enle}{\end{sl}\end{lemm}}
\newcommand{\bedef}{\begin{defi}\begin{sl}}
\newcommand{\eddef}{\end{sl}\end{defi}}
\newcommand{\bete}{\begin{teor}\begin{sl}}
\newcommand{\ente}{\end{sl}\end{teor}}
\newcommand{\beos}{\begin{osse}\begin{rm}}
\newcommand{\eddos}{\end{rm}\end{osse}}
\newcommand{\beas}{\begin{assu}\begin{rm}}
\newcommand{\eddas}{\end{rm}\end{assu}}
\newcommand{\bepr}{\begin{prop}\begin{sl}}
\newcommand{\empr}{\end{sl}\end{prop}}
\newcommand{\bepro}{\begin{prob}\begin{rm}}
\newcommand{\empro}{\end{rm}\end{prob}}
\newcommand{\bede}{\begin{defin}\begin{sl}}
\newcommand{\edde}{\end{sl}\end{defin}}
\newcommand{\beco}{\begin{coro}\begin{sl}}
\newcommand{\enco}{\end{sl}\end{coro}}

\newcommand{\quand}{\quad\text{and}\quad}

\newcommand{\quext}{\quad\text}

\newcommand{\RR}{\mathbb{R}}

\newcommand{\beeq}[1]{\begin{equation}\label{#1}}
\newcommand{\eddeq}{\end{equation}}

\newcommand{\beeqa}[1]{\begin{eqnarray}\label{#1}}
\newcommand{\eddeqa}{\end{eqnarray}}

\newcommand{\beal}[1]{\begin{align}\label{#1}}
\newcommand{\eddal}{\end{align}}

\newcommand{\bespl}[1]{\begin{split}\label{#1}}
\newcommand{\edspl}{\end{split}}

\newcommand{\bega}[1]{\begin{gather}\label{#1}}
\newcommand{\edga}{\end{gather}}

\newcommand{\beeqax}{\begin{eqnarray*}}
\newcommand{\eddeqax}{\end{eqnarray*}}

\def\qed{\ifmmode 
  \else \leavevmode\unskip\penalty9999 \hbox{}\nobreak\hfill
  \fi
  \quad\hbox{\hskip.5em\vrule width.4em height.6em depth.05em\hskip.1em}}
\def\endproofsym{\qed}
\renewenvironment{proof}[1][Proof]{\trivlist\item[\hskip\labelsep{\hskip0pt
    {\normalfont\scshape#1.}\hskip .321429\parindent}]\ignorespaces}
{\endproofsym\endtrivlist}
\def\endnobox{\def\endproofsym{}\end{proof}\def\endproofsym{\qed}}

\newcommand{\no}{\nonumber}

\newcommand{\beeqao}{\begin{eqnarray}\no}
\newcommand{\bealo}{\begin{align}\no}
\newcommand{\besplo}{\begin{split}\no}
\newcommand{\begao}{\begin{gather}\no}

\newcommand{\duav}[1]{\langle{#1}\rangle}

\newcommand{\perogni}{\forall\,}
\newcommand{\esiste}{\exists\,}

\newcommand{\io}{\int_\Omega}

\newcommand{\OO}{_{\Omega}}
\newcommand{\oo}{_{\Omega}}

\newcommand{\fhi}{\varphi}

\newcommand{\lhs}{left-hand side}
\newcommand{\rhs}{right-hand side}

\DeclareMathOperator{\deriv}{d}

\DeclareMathOperator{\sign}{sign}

\DeclareMathOperator{\loc}{loc}

\newcommand{\HUVp}{H^1(0,T;V')}
\newcommand{\CZH}{C^0([0,T];H)}
\newcommand{\CZV}{C^0([0,T];V)}

\newcommand{\LDH}{L^2(0,T;H)}
\newcommand{\LDV}{L^2(0,T;V)}

\newcommand{\LIH}{L^\infty(0,T;H)}

\newcommand{\LDHD}{L^2(0,T;H^2(\Omega))}

\let\TeXchi\chi
\def\chi{{\setbox0 \hbox{\mathsurround0pt
$\TeXchi$}\hbox{\raise\dp0 \copy0 }}}

\newcommand{\betaciapo}{\widehat{\beta}}

\newcommand{\calX}{{\mathcal X}}

\newcommand{\calA}{{\mathcal A}}

\newcommand{\calE}{{\mathcal E}}
\newcommand{\calD}{{\mathcal D}}

\newcommand{\calN}{{\mathcal N}}

\newcommand{\calS}{{\mathcal S}}

\newcommand{\eee}{_{\epsilon}}

\newcommand{\hgiu}{\underline{h}}
\newcommand{\sigmasu}{\overline{\sigma}}
\newcommand{\fhigiu}{\underline{\varphi}}

\newcommand{\dit}{\deriv\!t}
\newcommand{\dis}{\deriv\!s}
\newcommand{\dix}{\deriv\!x}

\newcommand{\ddt}{\frac{\deriv\!{}}{\dit}}

\newenvironment{bettirev}{\color{blue}}{\color{black}}
\newcommand{\bber}{\begin{bettirev}}
\newcommand{\eber}{\end{bettirev}}

\newenvironment{michelarev}{\color{red}}{\color{black}}
\newcommand{\III}{\begin{michelarev}}
\newcommand{\EEE}{\end{michelarev}}

\numberwithin{equation}{section}
\begin{document}

\title{On the long time behavior of a tumor growth model}

\author{%
Alain Miranville\\
Laboratoire de Math\'ematiques et Applications\\
UMR CNRS 7348, Equipe DACTIM-MIS\\
Universit\'e de Poitiers - SP2MI\\
Boulevard Marie et Pierre Curie\\
F-86962 Chasseneuil Futuroscope Cedex, France\\
Fudan University (Fudan Fellow), Shanghai, China\\
Xiamen University\\
School of Mathematical Sciences\\
Xiamen, Fujian, China\\
E-mail: {\tt Alain.Miranville@math.univ-poitiers.fr}\\
\and
Elisabetta Rocca\\
Dipartimento di Matematica, Universit\`a di Pavia, and IMATI - C.N.R.,\\
Via Ferrata~5, 27100 Pavia, Italy\\
E-mail: {\tt elisabetta.rocca@unipv.it}\\
\and
Giulio Schimperna\\
Dipartimento di Matematica, Universit\`a di Pavia, and IMATI - C.N.R.,\\
Via Ferrata~5, 27100 Pavia, Italy\\
E-mail: {\tt giusch04@unipv.it}
}


\maketitle
\begin{abstract}
 We consider the problem of the long time dynamics for a diffuse interface model for tumor growth.
 The model describes the growth of a tumor surrounded by host tissues in the presence of a nutrient 
 and consists in a Cahn-Hilliard-type equation for the tumor phase coupled with a reaction-diffusion 
 equation for the nutrient concentration. We prove that, under physically motivated assumptions
 on parameters and data, the corresponding initial-boundary value problem generates a dissipative
 dynamical system that admits the global attractor in a proper phase space.
\end{abstract}


\vspace{2mm}

\noindent {\bf AMS (MOS) subject clas\-si\-fi\-ca\-tion:}~~35D30; 35K57; 35Q92; 35B41; 37L30; 92C17

%
%

\section{Introduction}
\label{sec:intro}

One of the main examples of complex systems studied nowadays both in the biomedical and 
in the mathematical literature refers to tumor growth processes.
In particular, there has been a recent surge in the development of phase field models for tumor growth.  
These models aim to describe the evolution of a tumor mass surrounded by healthy tissues 
by taking into account biological mechanisms such as proliferation of cells 
via nutrient consumption, apoptosis, chemotaxis and active transport of specific chemical 
species. In particular, we will consider here a model that fits into the framework
of {\sl diffuse interface}\/ models for tumor growth. In this setting the evolution of
the tumor is described by means of an order parameter $\fhi$ that represents the 
local concentration of tumor cells; the interface between the tumor and healthy cells,
rather than being represented as a surface, is seen as a (narrow) layer separating
the regions where $\fhi=\pm1$, with $\fhi=1$ denoting the tumor phase 
and $\fhi= -1$ the healthy phase. Note that in the case of an incipient tumor, i.e., 
before the development of quiescent cells, the representation of the tumor
growth process is often given by a Cahn--Hilliard equation \cite{CH} for $\fhi$
coupled with a reaction-diffusion equation for the nutrient $\sigma$ 
(cf., e.g., \cite{CLLW09, GLSS16, HZO, HKNZ15}). We will consider here this type of 
situation; we just mention the fact that more sophisticated models may distinguish
between different tumor phases (e.g., proliferating and necrotic), or, treating
the cells as inertia-less fluids, include the effects of fluid flow into the evolution 
of the tumor, leading to (possibly multiphase) Cahn-Hilliard-Darcy systems
\cite{CL10,GLSS16,WLFC08}.

In this work, our main purpose is to consider the long time dynamics of 
a Cahn-Hilliard-reaction-diffusion tumor growth model recently introduced in \cite{GLSS16}. 
On the other hand, in comparison with \cite{GLSS16}, we have neglected here the effects 
of chemotaxis and active transport (a more complete model including these
effects may be the topic of a future work). Namely, we consider the following
PDE system:
\begin{align}\label{CH1}
  & \fhi_t - \Delta \mu = (P \sigma - A) h(\fhi),\\
 \label{CH2}
  & \mu = - \Delta \fhi + \psi'(\fhi),\\
 \label{nutr}
  & \sigma_t - \Delta \sigma = - C \sigma h(\fhi) + B (\sigma_s - \sigma),
\end{align}
settled in $\Omega\times (0,+\infty)$, $\Omega$ being a smooth domain of $\RR^3$,
and complemented with the Cauchy conditions and with no-flux
(i.e., homogeneous Neumann) boundary conditions for all
unknowns. As already mentioned, $\fhi$ represents the tumor phase concentration, 
$\sigma$ is the concentration of a nutrient for the tumor cells (such as oxygen or glucose),
and $\mu$~is the chemical potential of the ``phase transition'' from healthy
to tumor cells. The parameters
$P, A, B, C$ are assumed to be strictly positive constants, $\sigma_c \in (0,1)$,
and, in order to
ensure dissipativity, some compatibility conditions will be needed
(cf.~Assumption~\ref{ass:diss} below and Subsection~\ref{subsec:hom} for a detailed 
explanation of the effects of such conditions).
More specifically, in applications, $P$ denotes the proliferation rate, $A$ 
the apoptosis rate, $C$ the nutrient consumption rate, and $B$ the nutrient supply rate.
The term $P \sigma h(\fhi)$ models the proliferation of tumor cells which is proportional 
to the concentration of the nutrient, the term $A h(\fhi)$ describes the apoptosis of tumor cells, 
and $C \sigma h(\fhi) $ models the consumption of the nutrient by the tumor cells. 
The constant $\sigma_s$ denotes the nutrient concentration in a pre-existing vasculature,
and $B(\sigma_s - \sigma)$ models the supply of nutrient from the blood vessels 
if $\sigma_s > \sigma$ and the transport of nutrient away from the domain $\Omega$ if $\sigma_s< \sigma$. 
Moreover, $\psi'$~stands for the derivative of a double-well potential~$\psi$
and $h$ is a {smooth} proliferation function. A typical example of 
potential, meaningful in view of applications, has the expression
\begin{equation}\label{double}
  \psi_{reg}(r) = \frac 14 (r^2-1)^2,
  \quad r \in \RR,
\end{equation}
but we may observe that in our analysis we can allow for more general 
regular potentials having at least cubic and at most exponential growth at infinity. 
Hence, the polynomial potentials normally associated to the Cahn-Hilliard
energy are also admissible here. On the other hand, we may not consider
here the so-called {\sl singular potentials}, e.g.~of logarithmic type,
that are also popular in connection with Cahn-Hilliard-based models
(see, e.g., \cite{MAIMS}, cf.~also \cite{FLRS} for an application of logarithmic potentials to 
multiphase tumor growth models).

Let us now give, without any claim of completeness, a short overview of the 
recent mathematical literature on diffuse-interface tumor growth models. Modeling 
tumor growth dynamics has recently become a major issue in applied mathematics 
(see, e.g., \cite{Ciarletta,CL10,WLFC08}). Numerical simulations of diffuse 
interface models for tumor growth have been carried out in several papers 
(see, e.g., \cite[Ch.~8]{CL10}); nonetheless, a rigorous mathematical theory
of the related systems of PDEs is still at its beginning and many important 
problems are still open. We may quote \cite{CGH, CGRS1, CGRS2, DFRSS17, FGR, FLR,  GLDirichlet, GLNeumann}
as mathematical references for Cahn-Hilliard-type models and 
\cite{BosiaContiGrasselli14,GLDarcy,JiangWuZheng14,LowengrubTitiZhao13} 
for models also including a transport effect described by Darcy's law.

A further class of diffuse interface models that also include chemotaxis and transport effects 
has been subsequently introduced (cf.~\cite{GLNS, GLSS16}); moreover in some cases the sharp interface 
limits of such models have been investigated generally by using formal asymptotic methods.
Rigorous sharp interface limits have been however obtained in some special cases 
(see, e.g., the two recent works \cite{MR, RS}).

On the other hand, the problem of characterizing the long time behavior of solutions to tumor 
growth models is still in its infancy. Up to our knowledge, the only reference available to date
for Cahn-Hilliard-reaction-diffusion models is the work \cite{FGR}, where existence of the global attractor 
is proved in a phase space characterized by an a priori bound on the physical energy.
However, the model considered in \cite{FGR} has some notable differences with respect
to the present one (cf.~\cite{HZO} and see also \cite{HKNZ15,WZZ}). 
In particular, in \cite{FGR} the right-hand sides of \eqref{CH1} and \eqref{nutr} 
contain the chemical potential $\mu$ and this type of coupling implies that
a total energy balance can actually be proved. 

In this work, we prove the dissipativity of the system and the existence of a global attractor for the 
dynamical system generated by solutions of the initial-boundary value problem for \eqref{CH1}-\eqref{nutr}
taking values in the natural phase space which basically consists of the pairs $(\fhi,\sigma)$ 
having finite physical energy (cf.~\eqref{defi:X} below). The main mathematical
difficulty in the proof stands in establishing the dissipativity of the dynamical process, i.e., existence
of a uniformly absorbing set. Indeed, differently from standard Cahn-Hilliard models, 
here the spatial mean of $\fhi$ (i.e., the total mass of the tumor) is not conserved in time,
but the tumor may grow or shrink in a way that is essentially prescribed by the \rhs\ 
of~\eqref{CH1} which can be seen as a source of tumor mass. It is then clear that, if this
\rhs\ remains, say, positive for large values of $\fhi$, then the mass of $\fhi$ may grow
indefinitely and there can be no absorbing set. For this reason, dissipativity is 
only expected to hold under suitable compatibility conditions between the proliferation
function $h$ and the various coefficients $A,B,C,P, \sigma_s$. Roughly speaking these conditions
(which are thoroughly discussed below, see for instance Remark~\ref{su:hgiu}) prescribe
that, for large positive (negative) values of $\fhi$, the \rhs\ of \eqref{CH1}
must become negative (respectively, positive) in such a way that the tumor concentration
is forced to remain bounded in the $L^\infty$-norm uniformly for large values of the 
time variable. For this reason we need to assume in particular that, at least for $\fhi<<-1$,
$h(\fhi)$ stays {\sl strictly}\/ negative (and not equal to $0$ as was generally assumed
in former contributions); otherwise we cannot prove a uniform bound from below on $\fhi$.
We finally observe that, in view of our choice of no-flux boundary conditions, 
spatially homogeneous solutions exist. Their behavior is analyzed in Subsection~\ref{subsec:hom} 
by means of simple ODE techniques and in particular this gives further evidence of the 
fact that in absence of compatibility conditions on the coefficients, dissipativity of
the process may fail.

The paper is organized ed as follows: in the next section, we list our assumptions on the coefficients and data, 
state the problem in a precise form and present our main results. Then, the last section is devoted 
to the corresponding proofs and to a discussion on the mentioned compatibility conditions and on the 
behavior of spatially homogeneous solutions.

%
%

\section{Main results}
\label{sec:main}

We let $\Omega$ be a smooth bounded domain of $\RR^3$
with boundary~$\Gamma$. For simplicity, but with no loss of generality,
we assume $|\Omega|=1$. We set $H:=L^2(\Omega)$ and $V:=H^1(\Omega)$. We will use the same symbols
$H$ and $V$ for denoting vector valued functions (we may write, for instance,
$\nabla\fhi \in H$). The standard scalar product in $H$ will be noted as 
$(\cdot,\cdot)$. Since the immersion $V\subset H$ is continuous
and dense, identifying $H$ with its topological dual $H'$ through the above scalar product
we obtain the {\sl Hilbert triplet}\/ $(V,H,V')$. The duality pairing 
between a generic Banach space $X$ and its dual $X'$ will
be generally noted as $\duav{\cdot,\cdot}$. We let $R$ denote a weak form 
of the Laplace operator with
Neumann boundary conditions. Namely, we set
\begin{equation}\label{defiA}
   R:V \to V', \qquad
    \duav{Rv,z}:=\io \nabla v\cdot \nabla z \,\dix.
\end{equation}
For a generic function (or functional) $v$ defined over $\Omega$, we will
note its spatial mean value as
\begin{equation}\label{mean}
  v\OO := \frac1{|\Omega|} (v,1)
   = (v,1),
\end{equation}
the latter equality holding since $|\Omega|=1$. For, say, $v\in V'$,
the above holds replacing scalar products with duality pairings. We also recall
the Poincar\'e-Wirtinger inequality
\begin{equation}\label{powi}
   \| v - v\OO \| \le c_\Omega \| \nabla v \| \quad
    \perogni v \in V.
\end{equation}
Next, for any $\zeta\in V'$ we set 
\bega{defiV0}
  V_0':=\{\zeta\in V': \zeta\OO=0\},
   \qquad V_0:=V\cap V_0'.
\end{gather}
The above notation $V_0'$ is suggested just for the
sake of convenience; indeed, we mainly see $V_0$, $V_0'$ as 
(closed) subspaces of $V$, $V'$, inheriting their norms, 
rather than as a pair of spaces in duality.

Clearly, $R$ maps $V$ onto $V_0'$ 
and its restriction to $V_0$ is an isomorphism of $V_0$ 
onto $V_0'$. We denote by ${\cal N}:V_0'\to V_0$ 
the inverse of $R$, so that for any $u\in V$ and
$\zeta\in V_0'$ there holds
\beeq{relaN}
  \duav{R u,{\cal N}\zeta}
  =\duav{R{\cal N}\zeta,u}
    =\duav{\zeta,u}.
\end{equation}
We can now introduce a set of assumptions on the coefficients and data
that will be kept for the rest of the paper, noting that
some results will in fact require more specific conditions.
\beas\label{ass:gen}
 The coefficients are assumed to satisfy
 \begin{equation}\label{hp:const}
   P,A,B,C > 0, \qquad \sigma_c \in (0,1).
 \end{equation}
 The configuration potential $\psi$ lies in $C^{1,1}_{\loc}(\RR)$. Moreover its derivative
 is decomposed as a sum of a monotone part $\beta$ and a linear perturbation:
 \begin{equation}\label{hp:psi1}
   \psi'(r) = \beta(r) - \lambda r, \quad \lambda \ge 0,~~r \in \RR.
 \end{equation}
 The monotone part $\beta$ is normalized so that $\beta(0)=0$ and further 
 complies with the growth condition
 \begin{equation}\label{hp:beta1}
   \esiste c_\beta > 0:~~|\beta(r)| \le c_\beta ( 1 + \psi(r) )~~\perogni r\in\RR,
 \end{equation}
 which is more or less equivalent to asking $\psi$ to have
 at most an exponential growth at infinity. It will
 also be convenient to indicate by $\betaciapo$ the antiderivative of $\beta$ such
 that $\betaciapo(0)=0$. It then follows that $\betaciapo$ takes only
 nonnegative values; moreover, from \eqref{hp:psi1}, it turns out that
 $\psi(r) = \betaciapo(r) - \lambda r^2 / 2 + K$ for all $r\in \RR$,
 where $K$ is an integration constant which, thanks to \eqref{hp:beta1},
 can be chosen in such a way that $\min\psi = 0$.
 In order to avoid degenerate situations (such as $\beta=\psi\equiv0$, $\lambda=0$) 
 we also ask a minimal growth condition at infinity for $\psi$, i.e.~that
 \begin{equation}\label{hp:psi}
   \liminf_{|r|\nearrow \infty} \frac{\psi(r)}{|r|}
    =: \ell \in (0, +\infty]
 \end{equation}
 Only for the sake of proving uniqueness, condition \eqref{hp:beta1} has to be 
 slightly reinforced: we ask that there exists 
 $c>0$ such that 
 \begin{equation}\label{ass:beta1}
   | \beta(r) - \beta(s) | 
    \le c | r - s | \big( 1 + |\beta(r)| + |\beta(s)| \big)
    \quad \perogni r,s \in \RR.
 \end{equation}
 Note that this is still consistent with having at most an exponential 
 growth of $\beta$. 
 Next, we assume that $h$ is in $C^1(\RR)$, increasingly monotone 
 and it satisfies at least $h(-1)=0$ and $h(r)\equiv 1$ 
 for all $r\ge 1$. Moreover, we ask that there exist $\hgiu\ge 0$ 
 and $\fhigiu\le -1$ such that $h(r)\equiv -\hgiu$ for all $r\le \fhigiu$.
 Note that, as a consequence, $h$ is globally Lipschitz continuous. 
 Finally, we assume the initial data to satisfy
 \begin{align}\label{hp:init1}
   & \sigma_0 \in L^\infty(\Omega), \qquad
    0 \le \sigma_0 \le 1~~\text{a.e.~in }\,\Omega,\\
  \label{hp:init2}
   & \fhi_0 \in V, \qquad \psi(\fhi_0) \in L^1(\Omega).
 \end{align}
 We note that the second condition in \eqref{hp:init1} is not strictly
 necessary for proving existence. On the other hand, it makes sense to assume it
 in view of the physical interpretation of $\sigma$ as a nutrient concentration.
\eddas
\beos\label{su:h}
 %
 The simplest situation of a function $h$ satisfying the above assumption is given by
 the ``symmetric'' case corresponding to
 $\hgiu=0$ and $\fhigiu=-1$. On the other hand we will see in what follows that
 dissipativity of trajectories may not hold in such a case. This motivates
 our choice to consider the possibility of having $\hgiu>0$.
\eddos
\beos\label{su:h2}
 As mentioned in the introduction, it would also be significant to consider the case
 when $h(\fhi) = k \fhi + h_0(\fhi)$, where $k>0$ and $h_0$ is smooth and uniformly
 bounded; namely, $h$ is decomposed as a main linear part plus a bounded perturbation.
 This situation is somehow simpler because, at least as long as we can guarantee
 that $P\sigma - A < 0$, the linear part of $h$ drives some mass dissipation effect
 in \eqref{CH1}. 
\eddos
\beos\label{su:psi}
 As will be clear in a while when we discuss dissipativity, condition
 \eqref{hp:init2} corresponds to finiteness of the initial value
 of the ``physical'' energy (cf.~\eqref{energy} below).
 In particular,
 if $\psi$ grows at infinity as a polynomial of (possibly large)
 degree $p$, then the latter of \eqref{hp:init2} essentially 
 prescribes that $\fhi_0\in L^p(\Omega)$. 
\eddos
\beos\label{su:psi2}
 An explicit expression of a potential satisfying our hypotheses and having
 very slow (linear) growth at infinity is the following:
 \begin{equation}\label{ex:psi}
   \displaystyle \psi(r)=\begin{cases} \displaystyle
     \frac12 - r^2 \text{~~if }\,|r| \le \frac12,\\  \displaystyle
      (r-1)^2 \text{~~if }\,r \in\Big(\frac12,2\Big),\\  \displaystyle
      (r+1)^2 \text{~~if }\,r \in\Big(-2,-\frac12\Big),\\  \displaystyle
     2|r|-3\text{~~if }\,|r|\ge 2.
     \end{cases}
 \end{equation}
 Then, the conditions in Assumption~\ref{ass:gen} are satisfied with 
 $\lambda=\ell=2$. 
 On the other hand, we will see below that a potential like that in
 \eqref{ex:psi} is not suitable for having dissipativity, which seems
 to require a faster than cubic (but at most exponential) growth rate 
 at infinity. This growth rate is satisfied, for instance, by 
 the standard double-well potential 
 \eqref{double}.
\eddos
\noindent%
We are now ready to introduce our basic concept of weak solution:
\bede\label{def:weak}
 We say that a triplet $(\fhi,\mu,\sigma):(0,\infty)\times \Omega\to \RR^3$ 
 is a global weak solution to
 the tumor-growth model if the following conditions are satisfied:\\[1mm]
 (a)~for every $T>0$, there hold the regularity 
 properties
 \begin{align}\label{rego:fhi}
   & \fhi\in \HUVp \cap \CZV \cap \LDHD,\\
  \label{rego:beta}
   & \beta(\fhi)\in \LDH,\\
  \label{rego:mu}
   & \mu \in \LDV,\\
  \label{rego:sigma}
   & \sigma \in \HUVp \cap \CZH \cap \LDV \cap L^\infty(0,T;L^\infty(\Omega));
 \end{align}
 %
 %
 (b)~equations \eqref{CH1}-\eqref{nutr} are satisfied in the following weak sense:
 \begin{align}\label{CH1:w}
   & \fhi_t + R \mu = (P \sigma - A) h(\fhi)
   \quext{in }\,V',~~\text{a.e.~in }\,(0,\infty),\\
   \label{CH2:w}
   & \mu = R \fhi + \psi'(\fhi)
    \quext{in }\,H,~~\text{a.e.~in }\,(0,\infty),\\
   \label{nutr:w}
   & \sigma_t + R\sigma = - C \sigma h(\fhi) + B (\sigma_s - \sigma),
    \quext{in }\,V',~~\text{a.e.~in }\,(0,\infty);
  \end{align}
 (c)~there hold, a.e.~in~$\Omega$, the initial conditions
 \begin{equation}\label{init}
   \fhi|_{t=0}=\fhi_0, \qquad 
   \sigma|_{t=0}=\sigma_0.
 \end{equation}
\edde
\noindent%
Note that the homogeneous Neumann boundary conditions are now incorporated in the 
equations by definition of the operator~$R$ (cf.~\eqref{defiA}). Observe also that
\eqref{CH2:w} could in fact be interpreted as a pointwise relation (complemented
with an explicit boundary condition) thanks to the regularity \eqref{rego:fhi}.

\smallskip

Our first result is devoted to proving well-posedness in the class
of weak solutions:
\bete\label{teo:esi}
 Let\/ {\rm Assumption~\ref{ass:gen}} hold. Then the tumor-growth model admits 
 one and only one global in time weak solution 
 in the sense of\/ {\rm Definition~\ref{def:weak}}. Moreover, 
 for any $T>0$ there exists $\sigmasu_T \ge 1$ such that
 \begin{equation}\label{maxprinc}
   0 \le \sigma(t,x) \le \sigmasu_T, \quext{for a.e.~}\,(t,x) \in (0,T)\times \Omega,
 \end{equation}
 where we can take $\sigmasu_T$ independent of time if $B-C\hgiu >0$ and, in particular,
 $\sigmasu_T=1$ if $\hgiu=0$.
\ente
\noindent%
It is worth observing that existence and uniqueness hold without assuming any compatibility
conditions on the parameters $P, A, B, C, \sigma_c$. On the other hand,
as far as one wants to prove dissipativity of the dynamical process associated to
weak solutions, it seems necessary to take more restrictive assumptions. Note, for instance,
that \eqref{maxprinc} may allow the $L^\infty$-norm of $\sigma$ to increase in time. Hence,
we introduce a new
\beas\label{ass:diss}
 Let the parameters satisfy 
 \begin{align}\label{comp:1}
   & \hgiu > 0, \qquad B - C \hgiu > 0,\\
  \label{comp:1b}
   & \frac{B\sigma_s}{B - C \hgiu} < 1,\\
  \label{comp:2}
   & A - P \frac{B\sigma_s}{B - C\hgiu} > 0.
 \end{align}
 Let also $\beta$ have a superquadratic behavior at infinity, namely
 \begin{equation}\label{hp:psi4}
   \esiste \kappa_\beta > 0, C_\beta\ge 0, p_\beta>2:~~
    \beta(r)\sign r \ge \kappa_\beta |r|^{p_\beta} - C_\beta~~
     \perogni r \in \RR.   
 \end{equation}
\eddas 
\beos\label{su:hgiu}
 We remark that \eqref{comp:1}-\eqref{comp:1b} essentially prescribe
 $\hgiu$ to be {\sl strictly positive}, but {\sl small}. 
 The reason for such a condition will be clarified below. 
\eddos
\noindent
Our next result is actually devoted to proving that, if {\sl both}\/
Assumptions~\ref{ass:gen} and~\ref{ass:diss} hold, then weak solutions eventually
lie in a {\sl bounded absorbing set}\/ in a proper phase space. To define the latter, 
we introduce the usual Cahn-Hilliard energy functional
\begin{equation}\label{energy}
  \calE(\fhi)=\frac12 \| \nabla \fhi \|^2 + \io \psi(\fhi) \,\dix,
\end{equation}
arising as the sum of the interfacial and configurational energy. 
Then, we can define the ``energy space''
\begin{equation}\label{defi:X}
  \calX:=\big\{ (\fhi,\sigma) \in V \times L^\infty(\Omega):~\psi(\fhi) \in L^1(\Omega) \big\}
\end{equation}
and we correspondingly introduce the ``magnitude'' of an element $(\fhi,\sigma)\in\calX$ as
\begin{equation}\label{nor:X}
  \| (\fhi,\sigma) \|_{\calX} 
   := \| \fhi \|_V  + \| \sigma \|_{L^\infty(\Omega)} + \| \psi(\fhi) \|_{L^1(\Omega)}.
\end{equation}
Note that, in view of condition \eqref{maxprinc} (which holds with $\sigmasu$ independent
of $T$ thanks to \eqref{comp:1}), we already know that the component $\sigma$ of any 
weak solution stays bounded in $L^\infty(\Omega)$ uniformly in time. Observe also
that the quantity in \eqref{nor:X} is not a true norm due to the occurrence of the nonlinear
function $\psi$. On the other hand, convenience justifies the use of the
above notation.

\bigskip

We can state our second result about dissipativity of the dynamical process
generated by weak solutions:
\bete\label{teo:diss}
 Let\/ {\rm Assumptions~\ref{ass:gen} and \ref{ass:diss}} hold. Then there exists
 a positive constant $C_0$ independent of the initial data
 and a time $T_0$ depending only on the $\calX$-magnitude of the
 initial data such that any weak solution satisfies
 \begin{equation}\label{diss}
   \| (\fhi(t),\sigma(t)) \|_{\calX} \le C_0
    \quext{for every }\,t\ge T_0.
 \end{equation}
\ente
\noindent
Combining the above property with asymptotic compactness of trajectories,
we obtain the final result of this paper devoted to proving the 
existence of the global attractor. We refer the reader to, e.g., \cite{BV,MZ, Temam} for the related notions from the theory of 
infinite-dimensional dynamical systems.
\bete\label{teo:attr}
 Let\/ {\rm Assumptions~\ref{ass:gen}, \ref{ass:diss}} hold. Then the
 dynamical system generated by weak trajectories on the phase space
 $\calX$ admits the global attractor $\calA$. More precisely, $\calA$ 
 is a relatively compact subset of $\calX$ which is also bounded 
 in $H^2(\Omega) \times H^1(\Omega)$ and uniformly attracts the 
 trajectories emanating from any bounded set $B\subset \calX$.
\ente
\beos\label{reg:attr}
 In view of the fact that system \eqref{CH1}-\eqref{nutr} has a good
 parabolic structure, we expect the elements $(\fhi,\sigma)\in \calA$
 to be in fact smooth functions. More precisely their regularity may
 only be limited by the smoothness of the nonlinear functions $h$ and $\psi$. 
 In particular, if $h$ and $\psi$ are $C^\infty$, then the elements
 of the attractor are expected to be infinitely differentiable as well.
\eddos

%
%

\section{Proofs}
\label{sec:proofs}


\subsection{Proof of Theorem~\ref{teo:esi}: Well-posedness}
\label{subsec:exi}

{\bf A priori estimates.}~~%
The main ingredient of the proof of existence consists in a suitable set of
a priori estimates. To obtain them, we proceed here in a formal way by working
directly on equations \eqref{CH1}-\eqref{nutr}. The argument may 
however be easily justified within the framework of some regularization scheme
(e.g., Faedo-Galerkin). On the other hand, since the procedure
works similarly with related models (cf.~in particular \cite{GLNeumann}),
we just provide the highlights leaving details to the reader.

In what follows we will
note by $c>0$ and $\kappa>0$ some generic positive constants 
(whose specific value may vary on occurrence) depending only on 
the given parameters of the system (and neither on the initial data,
nor on any hypothetic approximation parameter).
The symbol $\kappa$ will be used in estimates from below.
Specific values of the constants will be noted as $c_i,\kappa_i$, 
$i\ge 1$. Constants depending on additional parameters will be noted 
using subscripts (e.g., $c_T$ if the constant depends on the
final time $T$).

To start with, we derive the basic boundedness properties for the nutrient. 
To this aim, we test \eqref{nutr} by $-\sigma_-$ (with $\sigma_-\ge 0$ denoting the 
{\sl negative part}\/ of $\sigma$) to deduce
\begin{equation}\label{co:01}
  \frac12 \ddt \| \sigma_- \|^2 
   + \| \nabla \sigma_- \|^2
   \le c \| \sigma_- \|^2.
\end{equation}
We used here the uniform boundedness of $h$ and the fact that
$B(\sigma_s-\sigma)$ is positive for $\sigma\le 0$ because $\sigma_s>0$.
Then, by \eqref{hp:init1} and the Gronwall lemma, we obtain
that $\sigma(t,x)\ge 0$ for (almost) every $t\ge 0$ and $x\in\Omega$.

To get an upper bound, we test \eqref{nutr} by $(\sigma-\sigmasu)_+$
with $\sigmasu\ge 1$ to be chosen below. Using the assumptions
on $h$ and performing standard manipulations, we deduce
\begin{align}\no
  & \frac12 \ddt \| (\sigma - \sigmasu)_+ \|^2 
   + \| \nabla ( \sigma - \sigmasu )_+ \|^2
   = - \io \big( ( B - C \hgiu ) \sigma  - B \sigma_s \big) ( \sigma - \sigmasu )_+ \, \dix\\
 \label{co:02}
 & \mbox{}~~~~~ 
   \le \io | B - C \hgiu | ( \sigma - \sigmasu )_+^2 \, \dix
   - \io \big( (B - C \hgiu )\sigmasu - B\sigma_s \big) ( \sigma - \sigmasu )_+ \, \dix\,.
\end{align}
We now have two cases. If $B - C \hgiu > 0$, then we can always choose $\sigmasu\ge 1$ 
large enough so that $(B - C \hgiu )\sigmasu - B\sigma_s \ge 0$. As a consequence, the
latter term on the \rhs\ is nonpositive and we can 
apply Gronwall's lemma to deduce that $\sigma(t,x)\le \sigmasu$ for 
a.e.~$(t,x)\in (0,\infty) \times \Omega$. Note that, if $\hgiu = 0$
the above certainly holds with $\sigmasu = 1$ in view of the 
fact that $\sigma_s<1$.

On the other hand, if $B - C \hgiu \le 0$, then the above procedure fails because 
we cannot control the last term in \eqref{co:02}. Nevertheless,
an $L^\infty$-estimate on $\sigma$ {\sl on finite times intervals}\/ can be obtained
also in that case. Indeed, one may test \eqref{nutr} by $\sigma^{p-1}$
(recall that we already know that in any case $\sigma\ge 0$) for 
a generic $p>1$. Then the boundedness of $h$ and easy computations give
\begin{equation}\label{co:03}
  \frac1p \ddt \| \sigma \|_{L^p(\Omega)}^p 
    \le c \big( 1 + \| \sigma \|_{L^p(\Omega)}^p \big)\,,
\end{equation}
with $c>0$ independent of $p$. Hence, setting $y_p:=\| \sigma \|_{L^p(\Omega)}^p$,
we obtain the differential inequality
\begin{equation}\label{co:04}
  (1+y_p)' \le cp (1+y_p), 
\end{equation}
whence 
\begin{equation}\label{co:05}
  \| \sigma(t) \|_{L^p(\Omega)}^p 
   \le 1 + y_p(t) 
   \le ( 1 + y_p(0) )e^{cpt}
   \le 2 e^{cpt}.
\end{equation}
Thus, taking the $1/p$-power and then letting $p\nearrow \infty$,
we get the desired conclusion. Summarizing, in any case we have obtained
\begin{equation}\label{st:11}
  \| \sigma \|_{L^\infty(0,T;L^\infty(\Omega))} \le c_T.
\end{equation}
This relation may be intended as an a priori estimate independent of any hypothetic
regularization parameter. Note that the constant on the \rhs\ is independent
of $T$ if $B - C \hgiu > 0$, and in particular it can be taken as $c_T=1$
if $\hgiu=0$.

\smallskip

As a next step, we derive the {\sl Energy estimate}\/ for the Cahn-Hilliard system.
This is the basic a priori information that any hypothetic weak solution 
is expected to satisfy. To obtain it,
we test \eqref{CH1} by $\mu$, \eqref{CH2} by $\fhi_t$ and sum up to obtain
\begin{equation}\label{en:11}
  \ddt \Big( \frac12 \| \nabla \fhi \|^2 + \io \psi(\fhi) \,\dix \Big)
   + \| \nabla \mu \|^2
   = \io (P\sigma - A) h(\fhi) \mu \,\dix.
\end{equation}
Then let us replace the expression for $\mu$ as given by \eqref{CH2}:
\begin{align}\no
  & \ddt \Big( \frac12 \| \nabla \fhi \|^2 + \io \psi(\fhi) \,\dix \Big)
   + \| \nabla \mu \|^2
   = \io (P\sigma-A) \big( h'(\fhi) | \nabla \fhi |^2 + h(\fhi) \beta(\fhi) \big) \,\dix\\ 
 \label{en:12}
  & \mbox{}~~~~~
   + \io \lambda (A - P\sigma) h(\fhi)\fhi \,\dix
   + P \io h(\fhi) \nabla \sigma \cdot \nabla \fhi\,\dix,
\end{align}
where $\psi'(\fhi)$ has been decomposed according to \eqref{hp:psi1}. Let us
now control the terms on the \rhs. First, as a
consequence of Assumption~\ref{ass:gen}, $|h(r)|+|h'(r)|\le c$
for every $r\in\RR$. Hence, using also \eqref{st:11},
\begin{equation}\label{co:20}
  \io (P\sigma-A) h'(\fhi) | \nabla \fhi |^2 \,\dix
   \le c \big( 1 + \| \sigma \|_{L^\infty(\Omega)} \big) \| \nabla \fhi \|^2
   \le c_T \| \nabla \fhi \|^2.
\end{equation}
Next, thanks to \eqref{hp:beta1},
\begin{equation}\label{co:21}
  \io (P\sigma-A) h(\fhi) \beta(\fhi) \,\dix 
   \le c \big( 1 + \| \sigma \|_{L^\infty(\Omega)} \big) 
   \bigg( 1 + \io \psi(\fhi) \,\dix \bigg)
    \le c_T + c_T \io \psi(\fhi) \,\dix.
\end{equation}
Finally, using also Young's inequality it is not difficult to deduce
\begin{equation}\label{co:22}
  \io \lambda (A - P\sigma) h(\fhi)\fhi \,\dix
   + P \io h(\fhi) \nabla \sigma \cdot \nabla \fhi\,\dix
  \le \frac12 \| \nabla \sigma \|^2 
  + c_T \big( 1 + \| \fhi \|_{L^1(\Omega)} + \| \nabla \fhi \|^2 \big).
\end{equation}
Note that the above constants $c_T$ depend on $T$ only through the 
$L^\infty$-norm of $\sigma$ (cf.~\eqref{st:11}).
In order to control the first term on the \rhs\ of
\eqref{co:22}, we test \eqref{nutr}
by $\sigma$. Then, straighforward calculations yield
\begin{equation}\label{co:23}
  \frac12 \ddt \| \sigma \|^2 
   + \| \nabla \sigma \|^2
   \le c \big( 1 + \| \sigma \|^2 \big).
\end{equation}
Summing \eqref{en:12} to \eqref{co:23} and using \eqref{co:20}-\eqref{co:22},
we arrive at 
\begin{align}\no
  & \ddt \Big( \frac12 \| \nabla \fhi \|^2 + \io \psi(\fhi) \,\dix + \frac12 \| \sigma \|^2 \Big)
   + \| \nabla \mu \|^2 + \frac12 \| \nabla \sigma \|^2\\
 \label{en:13}
  & \mbox{}~~~~~
  \le c_T \bigg( 1 + \| \fhi \|_{L^1(\Omega)} 
    + \| \nabla \fhi \|^2 + \| \sigma \|^2 
    + \io \psi(\fhi) \,\dix \bigg).
\end{align}
Now, using the growth assumption \eqref{hp:psi}, it is clear that, for 
some $\kappa,c>0$,  
\begin{equation}\label{psisu}
  \frac12 \| \nabla \fhi \|^2 + \io \psi(\fhi) \,\dix
   \ge \kappa \| \fhi \|_V - c \quad \perogni \fhi \in V.
\end{equation}
Hence, by Gronwall's lemma, \eqref{en:13} provides 
the following set of a priori estimates:
\begin{align}\label{st:12}
  &  \| \fhi \|_{L^\infty(0,T;V)} \le c_T,\\
  \label{st:13}
  &  \| \nabla \mu \|_{L^2(0,T;H)} \le c_T,\\
  \label{st:14}
  &  \| \psi(\fhi) \|_{L^\infty(0,T;L^1(\Omega))} \le c_T,\\
  \label{st:15}
  & \| \sigma \|_{L^2(0,T;V)\cap \LIH} \le c_T,
\end{align}
with $c_T$ as in \eqref{st:11}.

Next, integrating \eqref{CH2} over $\Omega$ and using once more 
\eqref{hp:beta1}, we deduce
\begin{equation}\label{co:31}
  | \mu\OO | = \bigg| \io \mu\,\dix \bigg| 
   = \bigg |\io ( \beta(\fhi) - \lambda \fhi )\,\dix \bigg|
    \le c_T \bigg( 1 + \io \psi(\fhi) \,\dix \bigg),
\end{equation}
where we have used \eqref{st:12} to control the $\lambda$-term. 
Recalling \eqref{st:14} we then infer
\begin{equation}\label{st:16}
  \| \mu\OO \|_{L^\infty(0,T)} \le c_T,
\end{equation}
which, combined with \eqref{st:13}, gives in turn
\begin{equation}\label{st:17}
  \| \mu \|_{L^2(0,T;V)} \le c_T.
\end{equation}
Now, testing \eqref{CH2} by $\beta(\fhi)$ and using \eqref{st:12}, \eqref{st:17}
and the monotonicity of $\beta$, it is a standard matter to deduce
\begin{equation}\label{st:18}
  \| \beta(\fhi) \|_{L^2(0,T;H)} \le c_T.
\end{equation}
Then, a comparison of terms in \eqref{CH2} and elliptic 
regularity results give
\begin{equation}\label{st:18b}
  \| \fhi \|_{L^2(0,T;H^2(\Omega))} \le c_T.
\end{equation}
Finally, we derive some estimates on the time derivatives of $\fhi$ and $\sigma$. 
Multiplying \eqref{CH1} by a generic nonzero test function $v\in V$ 
and using the previous estimates, we actually get
\begin{equation}\label{co:32}
  \duav{\fhi_t,v} = (\nabla \mu, \nabla v)
   + \io (P\sigma - A) h(\fhi) v\, \dix,
\end{equation}
whence estimates \eqref{st:12}, \eqref{st:17} and standard
manipulations yield
\begin{equation}\label{st:19}
  \| \fhi_t \|_{L^2(0,T;V')} \le c_T.
\end{equation}
Operating in an analogue way with equation \eqref{nutr} we similarly
obtain
\begin{equation}\label{st:19b}
  \| \sigma_t \|_{L^2(0,T;V')} \le c_T.
\end{equation}
\beos\label{add-reg}
 Using a more refined regularity argument in \eqref{CH2} and 3D
 Sobolev embeddings (see, e.g., \cite{MZ0}) one could
 improve \eqref{st:18}-\eqref{st:18b} up to
 \begin{equation}\label{st:18c}
   \| \beta(\fhi) \|_{L^2(0,T;L^6(\Omega))} 
    + \| \fhi \|_{L^2(0,T;W^{2,6}(\Omega))} \le c_T.
 \end{equation}
\eddos

\bigskip

\noindent%
{\bf Weak sequential stability.}~~%
We assume here to have a sequence of weak solutions $(\fhi_n,\mu_n,\sigma_n)$ satisfying
the a priori estimates obtained above uniformly with respect to the approximation 
parameter $n$. In other words, the constants $c$ or $c_T$ on the \rhs s 
of the bounds are assumed independent of $n$. We then prove that, up
to the extraction of subsequences, $(\fhi_n,\mu_n,\sigma_n)$ tends
in a suitable way to a triplet $(\fhi,\mu,\sigma)$ solving the tumor growth
model in the sense of Definition~\ref{def:weak} on the assigned but otherwise
arbitrary time interval $(0,T)$. This argument, generally noted as a ``weak
stability property'', may be seen as an abbreviated procedure for passing to
the limit in some approximation, for instance a Faedo-Galerkin scheme,
that may also involve the regularization of some terms (in particular of 
the function $\beta$). On the other hand, the procedure is so standard that we 
believe that giving very few highlights may suffice.

Actually, using the bounds \eqref{st:11}, \eqref{st:12}-\eqref{st:15}, 
\eqref{st:17}-\eqref{st:18b}, \eqref{st:19}-\eqref{st:19b}
and standard weak compactness argument, we are able to take a (nonrelabelled)
subsequence of $n$ such that $(\fhi_n,\mu_n,\sigma_n)\to (\fhi,\mu,\sigma)$
in the sense of weak or weak star convergence in proper Sobolev spaces. 
Moreover, using \eqref{st:19}, \eqref{st:19b}, and the Aubin-Lions lemma,
we obtain that $(\fhi_n,\sigma_n)$ tends to $(\fhi,\sigma)$ strongly in some
$L^p$-space, hence pointwise. This allows us to pass to the limit in the nonlinear
terms thanks to continuity of $h$ and $\beta$. In particular, we may observe
that, combining \eqref{st:18} with the pointwise convergence of $\fhi_n$
and using a generalized version of Lebesgue's dominated convergence theorem,
there follows
\begin{equation}\label{cons:18}
   \beta(\fhi_n) \to \beta(\fhi) 
    \quext{weakly in }\,L^2(0,T;H).
\end{equation}
Actually, even if in the approximation $\beta$ is replaced by some regularization $\beta_n$
the above property still works (with $\beta_n(\fhi_n)$ in place of $\beta(\fhi_n)$
on the \lhs) up to adaptations, provided that one assumes that $\beta_n$
tends to $\beta$ uniformly on compact subsets of $\RR$. 

\medskip

\noindent%
{\bf Uniqueness.}~~%
We give here a proof of uniqueness. A different (and somehow simpler) proof is given 
in \cite{GLR} (cf. also \cite{GLDirichlet}) in the case where $\psi$ has polynomial (of degree four) growth. On the other hand, the argument given
here works also for exponential $\psi$ (cf.~\eqref{ass:beta1}).
Assume to have two solutions $(\fhi_1,\mu_1,\sigma_1)$ and 
$(\fhi_2,\mu_2,\sigma_2)$ corresponding to two sets of initial data 
$(\fhi_{1,0},\sigma_{1,0})$ and $(\fhi_{2,0},\sigma_{2,0})$. 
Then the differences $(\fhi,\mu,\sigma):=(\fhi_1-\fhi_2,\mu_1-\mu_2,\sigma_1-\sigma_2)$ 
satisfy the following equations:
\begin{align}\label{CH1:d}
  & \fhi_t + R \mu = P \sigma h(\fhi_1)
    + ( P\sigma_2 - A )(h(\fhi_1)-h(\fhi_2)) 
  \quext{in }\,V',~~\text{a.e.~in }\,(0,\infty),\\
  \label{CH2:d}
  & \mu = R \fhi + \psi'(\fhi_1)-\psi'(\fhi_2)
   \quext{in }\,H,~~\text{a.e.~in }\,(0,\infty),\\
  \label{nutr:d}
  & \sigma_t + R \sigma = - C \sigma h(\fhi_1) -C\sigma_2(h(\fhi_1)-h(\fhi_2))- B \sigma,
   \quext{in }\,V',~~\text{a.e.~in }\,(0,\infty);
\end{align}
with  the initial conditions
\begin{equation}\label{init:d}
  \fhi|_{t=0}=\fhi_0, \qquad 
  \sigma|_{t=0}=\sigma_0,
\end{equation}
where $\fhi_0:=\fhi_{1,0}-\fhi_{2,0}$, $\sigma_0:=\sigma_{1,0}-\sigma_{2,0}$. 
In particular, integrating \eqref{CH1:d} over $\Omega$, we obtain
\begin{equation}\label{CH1:oo}
  (\fhi\OO)_t = \io P \sigma h(\fhi_1) \,\dix
    + \io ( P\sigma_2 - A )(h(\fhi_1)-h(\fhi_2)) \,\dix.
\end{equation}
Testing the above by $\fhi\OO$ and using the boundedness 
of $h$ and of $\sigma_2$ with the Lipschitz continuity of $h$,
we obtain
\begin{equation}\label{uni:11}
  \frac12 \ddt |\fhi\OO|^2
   \le c \big( | \fhi\OO |^2 + \| \sigma \|^2 + \| \fhi \|^2 \big).
\end{equation}
Next, let us take the difference of \eqref{CH1:d} and \eqref{CH1:oo}
and test it by $\calN(\fhi-\fhi\OO)$. Simpla calculations yield
\begin{equation}\label{uni:12}
  \frac12 \ddt \| \fhi - \fhi\OO\|_{V'}^2
   + \io \mu (\fhi - \fhi\OO) \, \dix
   \le c \big( \| \fhi - \fhi\OO\|_{V'}^2
    + \| \sigma \|^2 + \| \fhi \|^2 \big).
\end{equation}
Now, testing \eqref{CH2:d} by $\fhi-\fhi\OO$, we infer
\begin{align}\no
  \| \nabla \fhi \|^2  
   & = \io \mu (\fhi - \fhi\OO) \, \dix
   - \io ( \psi'(\fhi_1) - \psi'(\fhi_2) ) ( \fhi-\fhi\oo) \, \dix\\
 \label{uni:13}
   & \le \io \mu (\fhi - \fhi\OO) \, \dix
   + \fhi\oo \io ( \beta(\fhi_1) - \beta(\fhi_2) ) \, \dix
   + \lambda \| \fhi - \fhi\oo \|^2,
\end{align}
where we also used the decomposition \eqref{hp:psi1} and the monotonicity of 
$\beta$.

Next, testing \eqref{nutr:d} by $\sigma$, using the Lipschitz continuity of 
$h$ and performing standard manipulations, we deduce 
\begin{equation}\label{uni:14}
  \frac12 \ddt \| \sigma \|^2
   + \| \nabla \sigma \|^2  
   \le c \big( \| \sigma \|^2 + \| \fhi \|^2 \big).
\end{equation}
Combining \eqref{uni:11}-\eqref{uni:14}, we obtain
\begin{align}\no
  & \frac12 \ddt \big( | \fhi\oo|^2 
    + \| \fhi - \fhi\OO\|_{V'}^2
    + \| \sigma \|^2 \big)
    + \| \nabla \fhi \|^2  
    + \| \nabla \sigma \|^2 \\
 \label{uni:15}
   & \mbox{}~~~~~
   \le \fhi\oo \io ( \beta(\fhi_1) - \beta(\fhi_2) ) \, \dix
   + c \big( \| \fhi - \fhi\OO\|^2
    + \| \sigma \|^2 + | \fhi\oo |^2 \big).
\end{align}
In order to control the terms on the \rhs\ we first observe
that, thanks to the Poincar\'e-Wirtinger inequality and to Ehrling's lemma, 
\begin{equation}\label{uni:16}
  c \| \fhi - \fhi\OO\|^2
   \le \frac14 \| \nabla \fhi \|^2  
   + c \| \fhi - \fhi\OO\|^2_{V'}.
\end{equation}
To control the remaining term, we need to use 
assumption \eqref{ass:beta1} and then we derive
%
%
\begin{align}\no
  & \fhi\oo \io ( \beta(\fhi_1) - \beta(\fhi_2) ) \, \dix
    \le c | \fhi\oo | \io | \fhi | \big( 1 + |\beta(\fhi_1)| + |\beta(\fhi_2)| \big) \,\dix\\
 \no
   & \mbox{}~~~~~
   \le c | \fhi\oo | \| \fhi \| \big( 1 + \|\beta(\fhi_1)\| + \|\beta(\fhi_2)\| \big)\\
 \no
   & \mbox{}~~~~~
   \le c | \fhi\oo |^2  \big( 1 + \|\beta(\fhi_1)\|^2 + \|\beta(\fhi_2)\|^2 \big)
    + c \| \fhi  -\fhi\OO \|^2 + c | \fhi\OO |^2\\
 \label{uni:18}
   & \mbox{}~~~~~
   \le c | \fhi\oo |^2  \big( 1 + \|\beta(\fhi_1)\|^2 + \|\beta(\fhi_2)\|^2 \big)
    + c \| \fhi  -\fhi\OO \|_{V'}^2 + \frac14 \| \nabla \fhi \|^2.
\end{align}
Thanks to \eqref{uni:16} and \eqref{uni:18}, \eqref{uni:15} gives
\begin{align}\no
  & \frac12 \ddt \big( | \fhi\oo|^2 
    + \| \fhi - \fhi\OO\|_{V'}^2
    + \| \sigma \|^2 \big)
    + \frac12 \| \nabla \fhi \|^2  
    + \| \nabla \sigma \|^2 \\
 \label{uni:19}
   & \mbox{}~~~~~
   \le c | \fhi\oo |^2  \big( 1 + \|\beta(\fhi_1)\|^2 + \|\beta(\fhi_2)\|^2 \big)
   + c \big( \| \fhi - \fhi\OO\|_{V'}^2
    + \| \sigma \|^2 \big).
\end{align}
Then, using the regularity property \eqref{rego:beta} both for $\fhi_1$ and for $\fhi_2$
and applying Gronwall's lemma, we get uniqueness whenever
$(\fhi_{0,1},\sigma_{0,1})=(\fhi_{0,2},\sigma_{0,2})$. In the general case, we 
obtain the continuous dependence estimate 
\begin{align}\no
  & | (\fhi_1)\oo(t) - (\fhi_2)\oo(t) |^2 
    + \big\| \big( \fhi_1(t)- (\fhi_1)\oo(t) \big) - \big( \fhi_2(t)- (\fhi_2)\oo(t) \big) \big\|_{V'}^2
    + \| \sigma_1(t) - \sigma_2(t) \|^2 \\
 \label{cont:dep}
   & \mbox{}~~~~~
   \le C_T \Big( | (\fhi_{0,1})\oo - (\fhi_{0,2})\oo |^2 
    + \big\| \big(\fhi_{0,1}-(\fhi_{0,1})\oo\big) - \big(\fhi_{0,2}-(\fhi_{0,2})\oo\big) \big\|_{V'}^2
    + \| \sigma _{0,1} - \sigma_{0,2} \|^2 \Big),
\end{align}
for every $T>0$ and every $t\in(0,T]$, the constant $C_T>0$ depending on
the $\calX$-magnitude of the initial data and on $T$.


\subsection{Proof of Theorem~\ref{teo:diss}: Dissipativity}
\label{subsec:diss}

As a first step, we consider some auxiliary ODE's. Namely, we define $S_+$ 
and $S_-$ as the solutions to the following Cauchy
problems:
\begin{align}\label{eq:S+}
  & S_+' = - ( B - C \hgiu ) S_+ + B \sigma_s,\\
 \label{init:S+}
  & S_+(0) = 1,
\end{align}
and 
\begin{align}\label{eq:S-}
  & S_-' = ( - B - C ) S_- + B \sigma_s,\\
 \label{init:S-}
  & S_-(0) = 0.
\end{align}
Then we can readily compute
\begin{align}\label{S+}
  & S_+(t) = e^{-(B - C \hgiu )t} + \frac{B\sigma_s}{B-C\hgiu} \big( 1 - e^{-(B - C \hgiu )t} \big),\\
 \label{S-}
  & S_-(t) = \frac{B\sigma_s}{B+C} \big( 1 - e^{-( B + C )t} \big).
\end{align}
\bele\label{lemma1}
 Let the assumptions of\/ {\rm Theorem~\ref{teo:diss}} hold.
 Let $(\fhi,\sigma)$ be any weak solution to\/ \eqref{CH1}-\eqref{nutr}.
 Then we have
 \begin{equation}\label{subsuper}
   S_-(t) \le \sigma(t,x) \le S_+(t) 
    \quext{for every }\,t\ge 0~~\text{and }\,x\in \Omega.
 \end{equation}
\enle
\begin{proof}
We first recall that $\sigma(t,x)\ge 0$ for a.e.~$t\ge 0$, $x\in \Omega$
thanks to the minimum principle argument in the proof of Theorem~\ref{teo:esi}
(cf.~\eqref{co:01}).
Then, we can prove that $S_-$ is a subsolution, namely the first inequality
in \eqref{subsuper} holds. Taking the difference between \eqref{nutr}
and \eqref{eq:S-} we actually obtain
\begin{equation}\label{sub}
  (\sigma-S_-)' -\Delta (\sigma-S_-)
   = - B (\sigma - S_- ) - C (\sigma h(\fhi) - S_-),
\end{equation}
whence testing by $-(\sigma-S_-)_-$ and using the fact that $h\le 1$
we readily get the assert. Indeed, since $\sigma\ge 0$,
we notice that
\begin{equation}\label{contoS-}
  C (\sigma h(\fhi) - S_-) (\sigma-S_-)_-
  \le C (\sigma - S_-) (\sigma-S_-)_-
  \le 0.
\end{equation}
Analogously, the difference between \eqref{nutr}
and \eqref{eq:S+} gives
\begin{equation}\label{super}
  (\sigma-S_+)' -\Delta (\sigma-S_+) 
   = - B (\sigma - S_+ ) - C (\sigma h(\fhi) + S_+ \hgiu).
\end{equation}
Testing by $(\sigma-S_+)_+$, noting that
\begin{equation}\label{contoS+}
  - C (\sigma h(\fhi) + S_+ \hgiu ) (\sigma-S_+)_+
  \le - C (- \sigma \hgiu + S_+ \hgiu ) (\sigma-S_+)_+
  \le C \hgiu (\sigma-S_+)_+^2,
\end{equation}
and recalling \eqref{comp:1}, we easily obtain the second assertion.
\end{proof}
\noindent%
Recalling \eqref{comp:1b} and \eqref{comp:2}, we can take $\epsilon > 0$ to be a small number
satisfying
\begin{equation}\label{cond:epsi}
   2 \epsilon \le A - P \frac{B\sigma_s}{B - C\hgiu} 
    \quand \frac{B\sigma_s}{B - C\hgiu} + \frac{\epsilon}P < 1.
\end{equation}
We can then prove the following.
\bele\label{lemma2}
 Let the assumptions of\/ {\rm Theorem~\ref{teo:diss}} hold.
 Let $(\fhi,\sigma)$ be any weak solution in the sense of\/
 {\rm Definition~\ref{def:weak}.} Then there exist $T_1 > 0$ and 
 $C_1 > 0$ independent of the initial data such that
 \begin{align}\label{small:1}
   & \frac{B \sigma_s}{C+B} - \frac{\epsilon}P \le \sigma(t,x) \le \frac{B\sigma_s}{B - C\hgiu} + \frac{\epsilon}P 
    \quext{for all }\,t\ge T_1,~~\text{a.e.~}\,x \in \Omega,\\
  \label{small:2}
   & \| (\fhi(T_1),\sigma(T_1)) \|_{\calX} \le C_1 \big( 1 + \| (\fhi_0,\sigma_0) \|_{\calX} \big).
 \end{align}
\enle
\begin{proof}
Thanks to \eqref{subsuper} the component $\sigma$ evolves between the 
subsolution~$S_-$ and the supersolution $S_+$. Then, a simple computation
based on \eqref{S+}-\eqref{S-} shows that \eqref{small:1} holds provided
that we choose
\begin{equation}\label{T1}
  T_1:=\max\bigg\{ \frac{1}{B+C} \log\Big( \frac{B\sigma_sP}{\epsilon(B+C)} \Big), 
   \frac{1}{B-C\hgiu} \log\Big( \frac{P(B-C\hgiu-B\sigma_s)}{\epsilon(B-C\hgiu)} \Big) \bigg\}.
\end{equation}
Notice in particular that the argument of the second logarithm is strictly
positive thanks to assumption~\eqref{comp:1b}. Next, to prove \eqref{small:2},
it suffices to repeat the a priori estimates of Subsec.~\ref{subsec:exi}.
We may incidentally notice that the constant $c_T$ in \eqref{en:13} can now be taken
independent of $T$ thanks to Lemma~\ref{lemma1}. Anyway, integrating \eqref{en:13}
over the time interval $(0,T_1)$ and applying once more the Gronwall lemma,
we readily obtain the assertion.
\end{proof}
\noindent%
{\bf Proof of Theorem~\ref{teo:diss}.}~~
We start again from relation \eqref{en:12}, which we will now consider for
$t\ge T_1$. Hence, in particular we can take advantage of the second inequality 
in \eqref{small:1}. As a consequence, we can observe that,
thanks to \eqref{cond:epsi},
\begin{equation}\label{calSe2a}
  \sigma \le \frac{B\sigma_s}{B - C\hgiu} + \frac{\epsilon}P
   \Rightarrow P\sigma - A \le P \frac{B\sigma_s}{B - C\hgiu} + \epsilon - A
   \le - \epsilon.
\end{equation}
Consequently, for $t\ge T_1$ \eqref{en:12} implies
the following inequality:
\begin{align}\no
  & \ddt \Big( \frac12 \| \nabla \fhi \|^2 + \io \psi(\fhi) \,\dix \Big)
   + \| \nabla \mu \|^2
   + \epsilon \io h'(\fhi) | \nabla \fhi |^2 \, \dix\\
 \label{en:12d}
  & \mbox{}~~~~~
   + \io (A-P\sigma) h(\fhi) \beta(\fhi) \,\dix
   \le \io \lambda (A - P\sigma) h(\fhi)\fhi \,\dix
   + P \io h(\fhi) \nabla \sigma \cdot \nabla \fhi\,\dix.
\end{align}
Now, the terms on the \rhs\ can be controlled as in \eqref{co:22}.
On the other hand, using Assumptions~\ref{ass:gen} and \ref{ass:diss}
(and in particular the facts that $\hgiu$ is strictly positive and that $\beta(\fhi)$ has
the same sign as $\fhi$), it is not difficult to check that 
\begin{equation}\label{co:41z}
   h(\fhi) \beta(\fhi) \ge \kappa | \beta(\fhi) | - c,
\end{equation}
whence the latter term on the \lhs\ of \eqref{en:12d} gives
\begin{equation}\label{co:41}
   \io (A-P\sigma) h(\fhi) \beta(\fhi) \,\dix
   \ge \kappa\epsilon \| \beta(\fhi) \|_{L^1(\Omega)}
   - c,
\end{equation}
so that \eqref{en:12d} implies the differential inequality
\begin{equation}\label{en:12e}
  \ddt \Big( \frac12 \| \nabla \fhi \|^2 + \io \psi(\fhi) \,\dix \Big)
   + \| \nabla \mu \|^2
   + \kappa\epsilon \| \beta(\fhi) \|_{L^1(\Omega)}
  \le \frac12 \| \nabla \sigma \|^2 
   + c \big( 1 + \| \fhi \|_{L^1(\Omega)} + \| \nabla \fhi \|^2 \big).
\end{equation}
Adding \eqref{co:23} to the above relation, we arrive at
\begin{align}\no
  & \ddt \Big( \frac12 \| \nabla \fhi \|^2 + \io \psi(\fhi) \,\dix 
   + \frac12 \| \sigma \|^2 \Big) + \frac12 \| \nabla \sigma \|^2
   + \| \nabla \mu \|^2\\
  \label{en:12f}
  & \mbox{}~~~~~  + \kappa\epsilon \| \beta(\fhi) \|_{L^1(\Omega)}
    \le  c \big( 1 + \| \fhi \|_{L^1(\Omega)} + \| \nabla \fhi \|^2 \big),
\end{align}
where the norm of $\sigma$ on the \rhs\ of \eqref{co:23} has disappeared
because we now know that $0\le\sigma\le 1$ almost everywhere.

Next, let us multiply \eqref{CH2} by $-\Delta \fhi$. We deduce
\begin{equation}\label{en:31}
  \| \Delta \fhi \|^2 
   + \io \beta'(\fhi) | \nabla \fhi |^2 \,\dix 
  \le ( \nabla \fhi, \nabla \mu ) 
  - \lambda (\fhi, \Delta \fhi)
  \le ( \nabla \fhi, \nabla \mu ) 
   + \frac12 \| \Delta \fhi \|^2 
   + \frac{\lambda^2}2 \| \fhi \|^2.
\end{equation}
Correspondingly, testing \eqref{CH1} by $\fhi$ we obtain
\begin{equation}\label{en:32}
  \frac12 \ddt \| \fhi \|^2 
   + ( \nabla \fhi, \nabla \mu ) 
   + \io ( A - P \sigma ) h(\fhi) \fhi \, \dix = 0,
\end{equation}
whence in particular
\begin{equation}\label{en:33}
  \frac12 \ddt \| \fhi \|^2 
   + ( \nabla \fhi, \nabla \mu ) 
   \le c \big( 1 + \| \fhi \|^2 \big).
\end{equation}
Adding \eqref{en:31} and \eqref{en:33} to \eqref{en:12f}
and adding also the inequality $\frac12 \| \sigma \|^2 \le c$, 
neglecting some positive term on the \lhs, we obtain
\begin{align}\no
  & \ddt \Big( \frac12 \| \fhi \|_V^2 + \io \psi(\fhi) \,\dix 
   + \frac12 \| \sigma \|^2 \Big)
   + \frac{1}{2} \| \Delta \fhi \|^2 
   + \| \nabla \mu \|^2\\
  \label{en:12g}
  & \mbox{}~~~~~  + \kappa\epsilon \| \beta(\fhi) \|_{L^1(\Omega)}
   + \frac12 \| \sigma \|_V^2
   \le c \big( 1 + \| \fhi \|^2 + \| \nabla \fhi \|^2 \big).
\end{align}
Now, to control the \rhs, we first observe that 
\begin{equation}\label{en:33b}
  c \| \nabla \fhi \|^2 
   = c (- \Delta \fhi, \fhi )
   \le \frac14 \| \Delta \fhi \|^2 + c \| \fhi \|^2.
\end{equation}
Then, by virtue of assumption \eqref{hp:psi4},
for $\kappa,\epsilon$ as in \eqref{en:12g}, we have 
\begin{equation}\label{en:35}
  c \| \fhi \|^2 
   \le \frac{\kappa\epsilon}{2} \| \beta(\fhi) \|_{L^1(\Omega)}
   + c_{\kappa,\epsilon}.
\end{equation}
\beos\label{cresc:3}
 We point out that it may be possible to allow $p_\beta = 2$
 in \eqref{hp:psi4} at least in the case when $\kappa_\beta$ is large 
 enough. We leave the details to the reader. 
\eddos
\noindent%
Taking \eqref{en:33b} and \eqref{en:35} into account, \eqref{en:12g} gives
\begin{align}\no
 & \ddt \Big( \frac12 \| \fhi \|_V^2 + \io \psi(\fhi) \,\dix 
   + \frac12 \| \sigma \|^2 \Big)
   + \frac14 \| \Delta \fhi \|^2 
   + \frac{\kappa\epsilon}{2} \| \beta(\fhi) \|_{L^1(\Omega)}\\
  \label{en:12g2}
  & \mbox{}~~~~~  
   + \| \nabla \mu \|^2
   + \frac12 \| \sigma \|_V^2
   \le c.
\end{align}
Now, using \eqref{hp:psi4} again together with the continuous embedding
$H^2(\Omega)\subset L^\infty(\Omega)$, we notice that
\begin{equation}\label{en:36}
   \frac14 \| \Delta \fhi \|^2 
   + \frac{\kappa\epsilon}{2} \| \beta(\fhi) \|_{L^1(\Omega)}
   \ge \kappa \| \fhi \|_{H^2(\Omega)} - c
   \ge \kappa_1 \| \fhi \|_{L^\infty(\Omega)} - c.
\end{equation}
Let us then define 
\begin{equation}\label{defiZ}
  Z(r) := \betaciapo(r) + \betaciapo(-r), \quad \perogni r \ge 0,
\end{equation}
where $\betaciapo$ is the antiderivative of $\beta$ satisfying 
$\betaciapo(0)=0$ (hence in particular $\betaciapo$ is convex and
nonnegative due to Assumption~\ref{ass:gen}).
Noting that $Z$ is monotone over $[0,\infty§)$ with $Z(0)=0$, we have
\begin{equation}\label{orl:11}
  \io \betaciapo(\fhi) \, \dix
   \le \io Z(|\fhi|) \, \dix 
   \le \io Z\big( \| \fhi \|_{L^\infty(\Omega)} \big) \, \dix 
   = Z\big( \| \fhi \|_{L^\infty(\Omega)} \big).   
\end{equation}
As a consequence,
\begin{equation}\label{orl:12}
  \| \fhi \|_{L^\infty(\Omega)} 
   \ge Z^{-1} \bigg( \io \betaciapo(\fhi) \, \dix \bigg).
\end{equation}
Hence, recalling also \eqref{en:36}, relabelling
some constants, and rearranging some terms, \eqref{en:12g2} implies
\begin{align}\no
  & \ddt \Big[ \frac12 \| \fhi \|_V^2 
   + \frac12 \| \sigma \|^2 + \io \psi(\fhi) \,\dix \Big]
   + \frac{\kappa_3}2 \big( \| \fhi \|_V^2 + \| \sigma \|^2 \big)\\
 \label{en:12h}
  & \mbox{}~~~~~   
   + \kappa_1 Z^{-1} \bigg( \io \betaciapo(\fhi) \, \dix \bigg)
   + \kappa_2 \big( \| \Delta \fhi \|^2 + \| \nabla \mu \|^2 
   + \| \nabla \sigma \|^2 \big)
   \le c,
\end{align}
where the term $\frac{\kappa_3}2 \| \fhi \|_V^2$
has been added to both hands sides. Then its occurrence on the \rhs\ 
has been controlled essentially by repeating the procedure in 
\eqref{en:33b}-\eqref{en:35}. Now, for $K>0$ 
as in Assumption~\ref{ass:gen}, there holds 
\begin{equation}\label{orl:21}
  \betaciapo(r) 
   = \psi(r) + \frac\lambda2 r^2 - K
   \ge \psi(r) \quad \perogni |r|\ge \Big(\frac{2K}{\lambda}\Big)^{1/2}.
\end{equation}
As a consequence, for some $c>0$ we have
\begin{equation}\label{orl:22}
   \kappa_1 Z^{-1} \bigg( \io \betaciapo(\fhi) \, \dix \bigg)
    \ge  \kappa_1 Z^{-1} \bigg( \io \psi(\fhi) \, \dix \bigg) - c.
\end{equation}
Actually, to prove this relation it suffices to split the integration
domain $\Omega$ into the sets where $|\fhi|$ is smaller and respectively 
larger than $\left(\frac{2K}{\lambda}\right)^{1/2}$ and to use 
\eqref{orl:21}.

Thanks to the above relations, \eqref{en:12h} takes now the form
\begin{equation}\label{est:diss}
   \ddt \big( \calE_1 + \calE_2 \big) 
    + \kappa_3 \calE_1 + \kappa_1 Z^{-1}( \calE_2 )
    + \kappa_2 \calD \le c_1,
\end{equation}
where we have set
\begin{align}\label{defi:E12}
  & \calE_1:=  \frac12 \big( \| \fhi \|_V^2 + \| \sigma \|^2\big),
   \qquad \calE_2:= \io \psi(\fhi)\,\dix,\\
 \label{defi:D}
  & \calD := \| \Delta \fhi \|^2 + \| \nabla \mu \|^2 
   + \| \nabla \sigma \|^2
\end{align}
and we can notice that the above quantities are nonnegative. In
order to prove that the above differential inequality is dissipative,
we first observe that, as a consequence of \eqref{hp:beta1},
\begin{equation}\label{orl:31}
  \frac{|\beta(r)|}{\betaciapo(r)}\le c
   \quext{for sufficiently large }\,|r|. 
\end{equation}
whence, recalling \eqref{defiZ}, it is easy to deduce, for some $c\ge 0$,
\begin{equation}\label{orl:32}
  Z(r) \le c + e^{c r} \quad \perogni r \ge 0
\end{equation}
and, in turn, passing to inverse functions,
\begin{equation}\label{orl:33}
  Z^{-1}(r) \ge \kappa \ln ( y - c)  \quad \perogni r \ge \bar r,
\end{equation}
where $\bar r$ is some computable positive number. 
The above implies
\begin{equation}\label{orl:33b}
  \kappa_1 Z^{-1}(r) \ge \kappa_4 \ln ( y + 1) - c  \quad \perogni r \ge 0,
\end{equation}
so that inequality \eqref{est:diss} takes the form
\begin{equation}\label{est:diss:2}
  \ddt \big( \calE_1 + \calE_2 \big) 
   + \kappa_3 \calE_1 + \kappa_4 \ln (\calE_2 + 1) 
   \le c_2,
\end{equation}
and, using subadditivity of the logarithm,
\begin{equation}\label{est:diss:3}
  \ddt \big( \calE_1 + \calE_2 \big) 
   + \kappa_5 \ln (\calE_1 + \calE_2 + 1) 
   \le c_3,
\end{equation}
which is a dissipative differential inequality 
and implies the desired condition \eqref{diss}. Actually, 
it can be easily checked that there exists a 
finite and computable time $T_0\ge T_1$ depending only 
on the ``energy'' (in the sense of~\eqref{nor:X})
of the initial data such that for every $t\ge T_0$ there holds 
\begin{equation}\label{orl:24}
  \kappa_5 \ln (\calE_1 + \calE_2 + 1)  \le 2 c_3, \quext{i.e.\ }\,
   \calE_1 + \calE_2 \le e^{\frac{2c_3}{\kappa_5}} - 1.
\end{equation}
Indeed, if condition \eqref{orl:24} is violated, then the time derivative
of $\calE_1 + \calE_2$ is less than $-c_3$, implying that $\calE_1 + \calE_2$
decreases at least linearly with time until \eqref{orl:24} 
starts holding after some computable time $T_0$. Relation
\eqref{diss} is then an immediate consequence of \eqref{orl:24}.


\subsection{Proof of Theorem~\ref{teo:attr}: Attractor}
\label{subsec:attr}

Thanks to the dissipativity property of Theorem~\ref{teo:diss}, we only need
to show asymptotic compactness of solutions. To this
aim, we prove a further regularity estimate. As above, we will directly
work on system \eqref{CH1}-\eqref{nutr}, being intended that this formal
procedure may be justified within some approximation scheme. In what follows the 
various constants $c$ will be allowed to depend on the $\calX$-radius 
$C_0$ (cf.~\eqref{diss}) of the absorbing set.

That said, we first test \eqref{CH1} by $\mu_t$. Then, integrating by parts in 
time the term on the \rhs, we get
\begin{align}\no
  & (\mu_t,\fhi_t) + \frac12 \ddt \| \nabla \mu \|^2
   + \ddt \io (A - P\sigma) h(\fhi) \mu\,\dix 
   = \io \mu \big( (A - P\sigma) h(\fhi) \big)_t \dix \\
  \no
  & \mbox{}~~~~~ = \io (A - P\sigma) h'(\fhi)\fhi_t \mu\,\dix 
   - \io P \sigma_t h(\fhi) \mu\,\dix\\
  \label{co:61}
  & \mbox{}~~~~~ \le c \big( \| \fhi_t \| + \| \sigma_t \| \big) \| \mu \|
    \le \frac12 \| \fhi_t \|^2 + \frac12 \| \sigma_t \|^2 + c \| \mu \|^2.
\end{align}
We used here the boundedness of $h$ and $h'$, and the fact 
that $0\le \sigma \le 1$. These conditions will be repeatedly used again below
without further mentioning them. Next, we differentiate \eqref{CH2} in time 
and test the result by $\fhi_t$ to obtain
\begin{equation}\label{co:62}
  (\mu_t,\fhi_t) 
   = \| \nabla \fhi_t \|^2
   + \io \beta'(\fhi) \fhi_t^2 \,\dix
   - \lambda \| \fhi_t \|^2.
\end{equation}
Multiplying now \eqref{CH1} by $(1+2\lambda) \fhi_t$ we obtain
\begin{align}\no
  (1+2\lambda) \| \fhi_t \|^2
   & = - (1+2\lambda) (\nabla\mu,\nabla\fhi_t) 
   + (1+2\lambda) \io \big( P\sigma - A) h(\fhi) \fhi_t\,\dix\\
 \label{co:63}
   & \le \frac12 \| \nabla\fhi_t \|^2 + c_\lambda \| \nabla \mu \|^2
   + \lambda \| \fhi_t \|^2 + c_\lambda.
\end{align}
Finally, multiplying \eqref{nutr} by $2 \sigma_t$
and standardly controlling the \rhs, it is not difficult to deduce
\begin{equation}\label{co:64}
  \| \sigma_t \|^2
   + \ddt \| \nabla \sigma \|^2
   \le c.
\end{equation}
Taking the sum of relations \eqref{co:61}, \eqref{co:63}
and \eqref{co:64}, and using \eqref{co:62}, we arrive at
\begin{align}\no
  & \ddt \bigg[ \frac12 \| \nabla \mu \|^2
   + \| \nabla \sigma \|^2
   + \io (A - P\sigma) h(\fhi) \mu\,\dix \bigg]
  + \frac12\| \nabla \fhi_t \|^2
   + \io \beta'(\fhi) \fhi_t^2 \,\dix\\
  \label{co:65}
  & \mbox{}~~~~ + \frac12 \| \sigma_t \|^2
   + \frac12 \| \fhi_t \|^2
   \le c + c \| \nabla \mu \|^2
  + c \| \mu \|^2.
\end{align}
Now, using the Poincar\'e-Wirtinger inequality \eqref{powi}
we have
\begin{align}\no
  c \| \mu \|^2 
  & = c \| \mu - \mu\OO \|^2 + c \| \mu\OO \|^2 
  \le c \| \nabla \mu \|^2 + c \bigg| \io \psi'(\fhi) \, \dix \bigg|^2 \\
 \label{co:66} 
  & \le c \| \nabla \mu \|^2 + c + c \bigg| \io \psi(\fhi) \, \dix \bigg|^2
   \le c + c \| \nabla \mu \|^2,
\end{align}
where we have also used condition \eqref{hp:beta1} and the uniform bound
on the $L^1$-norm of $\psi(\fhi)$.

Then, noting as $\calE_3$ the sum of the terms in square brackets on the \lhs\ 
of \eqref{co:65}, we can observe that
\begin{align}\no
 \calE_3 & \ge \frac12 \| \nabla \mu \|^2
   + \| \nabla \sigma \|^2 - c \| \mu \|_{L^1(\Omega)}\\
 \no
 & \ge \frac12 \| \nabla \mu \|^2
   + \| \nabla \sigma \|^2 - c \| \mu - \mu\OO \|_{L^1(\Omega)}
   - c | \mu\OO |\\
 \no
 & \ge \frac12 \| \nabla \mu \|^2
   + \| \nabla \sigma \|^2 - c \| \nabla \mu \|
   - c - c \bigg| \io \psi(\fhi) \, \dix \bigg|\\
 \label{co:67}
 & \ge \frac14\| \nabla \mu \|^2
   + \| \nabla \sigma \|^2 - c_0,
\end{align}
where $c_0$ depends only on the uniform bound on the 
$\calX$-magnitude of the solution (cf.~\eqref{diss}) holding
for $t\ge T_0$.

Thanks to \eqref{co:66} and \eqref{co:67}, \eqref{co:65} gives
rise to the following inequality:
\begin{align}\no
  & \ddt \big( \calE_3 + c_0 )
  + \frac12\| \nabla \fhi_t \|^2
   + \io \beta'(\fhi) \fhi_t^2 \,\dix\\
  \label{co:65b}
  & \mbox{}~~~~ + \frac12 \| \sigma_t \|^2
   + \frac12 \| \fhi_t \|^2
   \le c + c \| \nabla \mu \|^2.
\end{align}
Now, coming back to \eqref{est:diss}, integrating it over the generic
time interval $(t,t+1)$, $t\ge T_1$, and recalling \eqref{defi:D},
we obtain
\begin{equation}\label{co:68}
  \int_t^{t+1} \big( \| \nabla \mu \|^2 
  + \| \nabla \sigma \|^2 \big) \,\dis
   \le c.
\end{equation}
Consequently, we can apply the uniform Gronwall lemma 
(see, e.g., \cite{Temam}) to \eqref{co:65b} to obtain
\begin{equation}\label{st:41}
  \| \mu(t) \|_V  
  + \| \sigma(t) \|_V
  \le C_1 \quad \perogni t \ge T_0 + 1,
\end{equation}
where $C_1>0$ is independent of the initial data. 
To get additional regularity on $\fhi$ it is then sufficient to go
back to \eqref{CH2} and apply standard elliptic regularity results
to obtain 
\begin{equation}\label{st:42}
  \| \fhi(t) \|_{H^2(\Omega)} \le C_2 \quad \perogni t \ge T_0 + 1,
\end{equation}
where $C_2>0$ is independent of the initial data. Properties \eqref{st:41}
and \eqref{st:42}, combined with the dissipativity proved in Theorem~\ref{teo:diss},
provide existence of the global attractor $\calA$ as well as its boundedness
in $H^2(\Omega)\times H^1(\Omega)$, which concludes the proof.


\subsection{Spatially homogeneous case}
\label{subsec:hom}

We give here some evidence of the fact that, if conditions \eqref{comp:1}-\eqref{comp:2}
do not hold, then dissipativity of the process may fail. To this aim
we will analyze the behavior of spatially homogeneous solutions. Indeed, in
view of the no-flux boundary conditions, these are particular solutions to
system \eqref{CH1}-\eqref{nutr} starting from spatially homogeneous initial data.
Then, let us denote by $X=X(t)$ and by $S=S(t)$ the 
spatially homogeneous versions of $\fhi$ and $\sigma$, respectively. 
In this setting, our problem reduces to the following ODE system 
for the vector variable $(X,S)$:
\begin{align}\label{CH:spat}
  & X' + (A - P S) h(X) = 0,\\
 \label{nutr:spat}
  & S' + C S h(X) + B (S - \sigma_s) = 0.
\end{align}
We can first observe that, if $\hgiu = 0$ and $X(0)<-1$,
equation \eqref{CH:spat} prescribes $X(t)$ to be conserved in time. 
In particular, there is no hope to prove that $X(t)$ eventually
lies in some bounded absorbing set. Indeed, if $X$ is large negative
at the initial time, then it remains like that forever.

Let us now assume $\hgiu>0$. Then we observe that
\begin{equation}\label{co:11a}
  B\sigma_s - ( C + B ) S \le S' \le B \sigma_s - ( B - C \hgiu) S.
\end{equation}
The first inequality implies that 
\begin{equation}\label{co:12a}
  S < \frac{B \sigma_s}{C+B} \Rightarrow S' > 0.
\end{equation}
For what concerns the second inequality, we have two cases. Let us 
first consider the situation when $C \hgiu \ge B$,
i.e.~\eqref{comp:1} does not hold. Let also the initial 
data be chosen in such a way that $X(0)<<0$ and $S(0)>>0$
(in such a way that $PS-A>0$). Then it follows
\begin{align}\label{CH:spat1}
  & X' = - (P S - A) \hgiu < 0,\\
 \label{nutr:spat1}
  & S' = B \sigma_s + (C\hgiu - B)S > 0
\end{align}
and both $|X|$ and $S$ go increasing forever. Note that in this situation,
even if we restrict ourselves  to the ``physical'' case
$S(0)\in [0,1]$, if $X(0)<-1$ then due to \eqref{nutr:spat1}
$S(t)$ eventually becomes larger than $1$; hence, the physical constraint
$S(t)\in[0,1]$ is not respected. 

In view of the above discussion, it looks reasonable to assume 
$\hgiu>0$ and \eqref{comp:1}. Under these conditions, the second inequality 
in \eqref{co:11a} implies
\begin{equation}\label{co:12b}
  S > \frac{B\sigma_s}{B - C\hgiu} \Rightarrow S' < 0.
\end{equation}
We can then define the region 
\begin{equation}\label{co:13}
   \calS:= \Big\{ (X,S) \in \RR^2:~ \frac{B \sigma_s}{C+B} \le S \le \frac{B\sigma_s}{B-C\hgiu} \Big\}
\end{equation}
and it follows from \eqref{co:12a} and \eqref{co:12b} that $\calS$ is positively 
invariant for the dynamical process generated by \eqref{CH:spat}-\eqref{nutr:spat}.
Now, if we want to keep the physical constraint $S(t)\in[0,1]$,
we need to assume $\frac{B\sigma_s}{B-C\hgiu}<1$, i.e.~\eqref{comp:1b}
(otherwise basically our results still hold provided that we allow
$S$ to take also values larger than $1$). In such a situation, we 
need to emphasize the role of \eqref{comp:2}. To this purpose, let us
assume that $X(0)>1$, which also implies $h(X)=1$. Then, \eqref{CH:spat}
reduces to 
\begin{equation}\label{CH:spat2}
  X' = (P S - A)
\end{equation}
and in this sense condition \eqref{comp:2} (which can be rewritten as
$\frac{A}{P}>\frac{B\sigma_s}{B-C\hgiu}$) prescribes 
that (if we reason in the $(X,S)$-plane with $X$ represented in
the horizontal axis), in the intersection between $\calS$ and 
the semiplane $\{X>1\}$, $X'$ stays negative (hence arbitrary growth
of $X$ is prevented, because trajectories tend to eventually
enter the region $\calS$). 

On the other hand, we can see that, when $\frac{A}{P}\le\frac{B\sigma_s}{B+C}$, dissipativity
cannot hold. Indeed if $S(0)\in \left[\frac{B \sigma_s}{C+B},\frac{B\sigma_s}{B-C\hgiu}\right]$
and $X(0)\ge 1$, then $X(t)$ is forced to increase forever, because 
$(X,S)$ can never leave the positively invariant region $\calS$
where, now, $X'>0$. On the other hand, the situation when 
$\frac{A}{P}\in \left(\frac{B \sigma_s}{C+B},\frac{B\sigma_s}{B-C\hgiu}\right]$ is
unclear, in the sense that, when $X>1$, in the ``upper'' part of the strip $\calS$,
$X'$ is positive, whereas $X'$ is negative in the ``lower'' part of $\calS$, so 
the evolution of $(X,S)$ may be more difficult to capture. Of
course, the behavior may be even more complicated once one 
considers general (i.e., not necessarily spatially homogeneous)
solutions to \eqref{CH1}-\eqref{nutr}, because in that case also 
equation \eqref{CH2} plays an important role
(whereas \eqref{CH2} ``disappears'' in the  spatially
homogeneous setting).

\bigskip
\noindent
{\bf Acknowledgments.}~~
This research has been performed in the framework of the project Fondazione Cariplo-Regione Lombardia  MEGAsTAR
``Matema\-tica d'Eccellenza in biologia ed ingegneria come acceleratore
di una nuova strateGia per l'ATtRattivit\`a dell'ateneo pavese''. The present paper
also benefits from the support of the MIUR-PRIN Grant 2015PA5MP7 ``Calculus of Variations'' for GS,
and of the GNAMPA (Gruppo Nazionale per l'Analisi Matematica, la Probabilit\`a e le loro Applicazioni)
of INdAM (Istituto Nazionale di Alta Matematica) for ER and~GS.




\begin{thebibliography}{99}


\bibitem{Ciarletta}
A.~Agosti, P.F.~Antonietti, P.~Ciarletta, M.~Grasselli, M.~Verani, 
{\sl A Cahn-Hilliard-type equation with application to tumor growth dynamics}, 
Math.\ Methods Appl.\ Sci., 
{\bf 40} (2017), 7598--7626.

\bibitem{BV} 
{A.V. Babin, M.I. Vishik}, {Attractors of evolution
equations}, {North-Holland}, {Amsterdam}, {1992}.



\bibitem{BosiaContiGrasselli14}
S.~Bosia, M.~Conti, M.~Grasselli,
{\sl On the {C}ahn--{H}illiard--{B}rinkman system},
{Commun.\ Math.\ Sci.},
{\bf 13},
(2015),
{1541--1567}. 

\bibitem{CH}
 J.W.~Cahn and J.E.~Hilliard, 
 {\sl Free energy of a nonuniform system. I.~Interfacial free energy},
 J.~Chem.\ Phys., 
 {\bf 28} (1958),
 258--267.

\bibitem{CGH}
P.~Colli, G.~Gilardi, D.~Hilhorst, 
{\sl On a Cahn--Hilliard type phase field model related to tumor growth}, 
{Discrete Contin.\ Dyn.\ Syst.}, 
{\bf 35} (2015), 2423--2442.

\bibitem{CGRS1}
P.~Colli, G.~Gilardi, E.~Rocca, J.~Sprekels, 
{\sl Vanishing viscosities and error estimate for a Cahn--Hilliard type phase-field system related to tumor growth}, 
{Nonlinear Anal.\ Real World Appl.}, {\bf 26} (2015), 93--108.

\bibitem{CGRS2}
P.~Colli, G.~Gilardi, E.~Rocca, J.~Sprekels, 
{\sl Asymptotic analyses and error estimates for a Cahn--Hilliard type phase field system modelling tumor growth}, 
{Discrete Contin.\ Dyn.\ Syst.\ Ser.~S.}, 
{\bf 10} (2017), 37--54.
 
\bibitem{CL10} V.~Cristini, J.~Lowengrub, 
Multiscale modeling of cancer. An integrated experimental and mathematical modeling approach,
Cambridge Univ.~Press, 2010. 

\bibitem{CLLW09}
V.~Cristini, X.~Li, J.S.~Lowengrub, S.M.~Wise, 
{\sl Nonlinear simulations of solid tumor growth using a mixture model: {invasion} and branching},
{J.~Math.\ Biol.}, 
{\bf 58} (2009), 723--763. 
 
\bibitem{DFRSS17}
M.~Dai, E.~Feireisl, E.~Rocca, G.~Schimperna, M.~Schonbek, 
{\sl Analysis of a diffuse interface model for multispecies tumor growth}, 
{Nonlinearity}, 
{\bf 30} (2017), 1639--1658.

\bibitem{FGR}
S.~Frigeri, M.~Grasselli, E.~Rocca, 
{\sl On a diffuse interface model of tumor growth}, 
{European J.~Appl.\ Math.}, 
{\bf 26} (2015), 215--243.

\bibitem{FLR}
S.~Frigeri, K.F.~Lam, E.~Rocca, 
{\sl On a diffuse interface model for tumour growth 
with non-local interactions and degenerate mobilities}. 
In: P.~Colli, A.~Favini, E.~Rocca, G.~Schimperna, J.~Sprekels~(eds.), Solvability,
Regularity, Optimal Control of Boundary Value Problems for PDEs,
pp.~217--254, Springer INdAM Series, Springer, Milan, 2017.

\bibitem{FLRS}
S.~Frigeri, K.F.~Lam, E.~Rocca, G.~Schimperna,
{\sl On a multi-species Cahn-Hilliard-Darcy tumor growth model with singular potentials},
Commun.\ Math.\ Sci., 
{\bf 16} (2018), 821--856.


\bibitem{GLDarcy}
H.~Garcke, K.F.~Lam, 
{\sl Global weak solutions and asymptotic limits of a Cahn--Hilliard--Darcy system modelling tumour growth}, 
{AIMS Mathematics}, {\bf 1} (2016), 
318--360.

\bibitem{GLDirichlet}
H.~Garcke, K.F.~Lam, 
{\sl Analysis of a Cahn--Hilliard system with non zero Dirichlet conditions modelling tumour growth with chemotaxis},
{Discrete Contin.\ Dyn.\ Syst.}, 
{\bf 37} (2017), 4277--4308.

\bibitem{GLNeumann}
H.~Garcke, K.F.~Lam, 
{\sl Well-posedness of a Cahn--Hilliard system modelling tumour growth with chemotaxis and active transport}, 
{European J.~Appl.\ Math.}, 
{\bf 28} (2017), 284--316.

\bibitem{GLNS} 
H.~Garcke, K.F.~Lam, {R.~N\"urnberg}, E.~Sitka, 
{\sl A multiphase Cahn--Hilliard--Darcy model for tumour growth with necrosis}, 
Math.\ Models Methods Appl.\ Sci., 
{\bf 28} (2018), 525--577.

\bibitem{GLR}
H.~Garcke, K.F.~Lam, E.~Rocca, {\it Optimal control of treatment time in a diffuse interface model for tumour growth},
Appl Math Optim, {\bf 78} (2018), 495--544.

\bibitem{GLSS16}
H.~Garcke, K.F.~Lam, E.~Sitka, V.~Styles, 
{\sl A {Cahn--Hilliard--Darcy} model for
tumour growth with chemotaxis and active transport}, 
Math.\ Models Methods Appl.\ Sci.,
{\bf 26} (2016),
1095--1148. 
 
 
\bibitem{HZO} 
A.~Hawkins-Daarud, K.G.~van der Zee, J.T.~Oden,
{\sl Numerical simulation of a thermodynamically consistent four-species tumor growth model}, 
Int.\ J.~Numer.\ Meth.\ Biomed.\ Engng., 
\textbf{28} (2011), 3--24.

\bibitem{HPZO} 
A.~Hawkins-Daarud, S.~Prudhomme, K.G.~van der Zee, J.T.~Oden,
{\sl Bayesian calibration, validation, and uncertainty 
quantification of diffuse interface models of tumor growth}, 
J.~Math.\ Biol., 
\textbf{67} (2013), 
1457--1485.

\bibitem{HKNZ15} 
D.~Hilhorst, J.~Kampmann, T.N.~Nguyen, K.G.~van der Zee, 
{\sl Formal asymptotic limit of a diffuse-interface tumor-growth model}, 
{Math.\ Models Methods Appl.\ Sci.}, 
{\bf 25} (2015), 1011--1043. 

\bibitem{JiangWuZheng14}
J.~Jiang, H.~Wu, S.~Zheng,
{\sl Well-posedness and long-time behavior of a non-autonomous {C}ahn--{H}illiard--{D}arcy 
system with mass source modeling tumor growth},
{J.~Differential Equations},
{\bf 259} (2015), 
{3032--3077}.

\bibitem{LowengrubTitiZhao13}
{J.S.~Lowengrub, E.~Titi, K.Zhao},
{\sl Analysis of a mixture model of tumor growth},
{European J.~Appl.\ Math.},
{\bf 24} (2013),
{691--734}. 

\bibitem{MR}
S.~Melchionna, E.~Rocca, 
{\sl Varifold solutions of a sharp interface limit of a diffuse interface model for tumor growth}, 
{Interfaces Free Bound.},
{\bf 19} (2018), 571--590. 

\bibitem{MAIMS}
A.~Miranville,
{\sl The Cahn-Hilliard equation and some of its variants},
{AIMS Mathematics},
{\bf 2} (2017), 479--544.

\bibitem{MZ0}
A.~Miranville, S.~Zelik,
{\sl Robust exponential attractors for Cahn-Hilliard type equations with singular potentials},
Math.\ Methods Appl.\ Sci.,
{\bf 27} (2004), 
545--582.

\bibitem{MZ}
A. Miranville, S. Zelik,
     {Attractors for dissipative partial differential equations in bounded and unbounded domains},
      in ``Handbook of Differential Equations: Evolutionary Equations, Vol. IV" (eds. C.M. Dafermos and M. Pokorny),
                Elsevier/North-Holland, 103--200, 2008.


\bibitem{RS}
E.~Rocca, R.~Scala,
{\sl A rigorous sharp interface limit of a diffuse interface model related to tumor growth},
{J.~Nonlinear Sci.}, 
{\bf 27} (2017), 847--872. 

\bibitem{Temam} R.~Temam, 
Infinite Dimensional Dynamical Systems in Mechanics and Physics, 
Springer, New York, 1997.

\bibitem{WLFC08} S.M.~Wise, J.S.~Lowengrub, H.B.~Frieboes, V.~Cristini, 
{\sl Three-dimensional multispecies nonlinear tumor growth--I: 
model and numerical method}, 
{J.~Theoret.\ Biol.}, 
{\bf 253} (2008), 524--543.

\bibitem{WZZ} 
X.~Wu, G.J.~van Zwieten, K.G.~van der Zee,
{\sl Stabilized second-order convex splitting schemes for Cahn-Hilliard
models with applications to diffuse-interface tumor-growth models},
Int.~J.~Numer.\ Meth.\ Biomed.\ Engng., 
\textbf{30} (2014), 180--203.




\end{thebibliography}
\end{document}